\documentclass[11pt]{article}
\usepackage{latexsym, graphics, amsfonts}

\newtheorem{thm}{Theorem}[section]
\newtheorem{cor}[thm]{Corollary}
\newtheorem{lem}[thm]{Lemma}
\newtheorem{prop}[thm]{Proposition}
\newtheorem{rem}{Remark}
\newcommand{\match}[1]{\mathcal{M}(#1)}

\newcommand{\Pf}{\mathrm{Pf}}

\newenvironment{proof}{{\bf Proof:}}{\hfill\rule{2mm}{2mm}}
\newcommand{\sgn}{\mathrm{sgn}}
\begin{document}

\title{Graphical Condensation Generalizations Involving Pfaffians and Determinants}%
\author{Eric H. Kuo\footnote{Department of Mathematical Sciences, George Mason University, Faifax, VA 22030.  E-mail:{\tt ekuo@gmu.edu.}}}%
\maketitle


\begin{abstract}
Graphical condensation is a technique used to prove combinatorial identities among numbers of perfect matchings of plane graphs. Propp and Kuo first applied this technique to prove identities for bipartite graphs. Yan, Yeh, and Zhang later applied graphical condensation to nonbipartite graphs to prove more complex identities. Here we generalize some of the identities of Yan, Yeh, and Zhang. We also describe the latest generalization of graphical condensation in which the number of perfect matchings of a plane graph is expressed as a Pfaffian or a determinant where the entries are also numbers of perfect matchings of subgraphs.
\end{abstract}
\maketitle
\section{Introduction}
In this article we will let $G=(V(G),E(G))$ be a simple graph. The set of perfect matchings will be denoted $\match{G}$, and the number of perfect matchings is $M(G) = |\match{G}|$. If $G$ is a weighted graph, the weight of a perfect matching is the product of the weights of the edges in the matching. The sum of the weights of perfect matchings of $G$ is also denoted $M(G)$. If $U$ is a subset of vertices in $V(G)$, then $G-U$ is the subgraph of $G$ induced by the vertices of $V\backslash U$. If $a$ is a vertex in $G$, then $G-a = G-\{a\}$. Finally, the length of a path (or cycle) is the number of edges in the path. A path is even or odd based on its length.

The following proposition first appeared in an email sent to the ``Domino Forum":
\begin{prop}\label{prop:int1}
Let $G$ be a plane graph with four vertices $a,b,c,d$ that appear in that cyclic order on a face of $G$. Then
\begin{eqnarray}
\lefteqn{M(G) M(G-\{a,b,c,d\}) + M(G-\{a,c\}) M(G-\{b,d\}) =} \nonumber \\
& & M(G-\{a,b\}) M(G-\{c,d\}) + M(G-\{a,d\}) M(G-\{b,c\}).
\end{eqnarray}
\end{prop}
This proposition generalizes results first reported by Propp~\cite{Propp} and Kuo~\cite{Kuo} for which $G$ was a bipartite graph. A combinatorial proof of Proposition~\ref{prop:int1} was first published by Yan, Yeh, and Zhang~\cite{YYZ}. They first proved the following generalization from which Proposition~\ref{prop:int1} follows as a corollary:
\begin{thm}\label{thm:yyz}
Let $G$ be a plane weighted graph with $2n$ vertices. Let $2k$ vertices $a_1,b_1,a_2,b_2, \ldots, a_k,b_k$ (where $2\leq k \leq n$) appear in cyclic order on a face of $G$. Let $A=\{a_1,\ldots, a_k\}$, and $B=\{b_1,\ldots, b_k\}$. Then for any $j=1,2,\ldots, k$, we have 
\begin{equation}\label{eqn:yyz}
\sum_{W\subseteq B, |W| \:\mathrm{even}} M(G-W)M(G-A-\overline{W}) 
= \sum_{Y\subseteq B, |Y| \:\mathrm{odd}} M(G-a_j-Y)M(G-A\backslash\{a_j\}-\overline{Y})
\end{equation}
where $\overline{W}=B\backslash W$ and $\overline{Y}=B\backslash Y$, the first sum ranges over all even subsets $W$ of $B$, and the second sum ranges over all odd subsets $Y$ of $B$.
\end{thm}

This theorem was proved using a lemma from Ciucu~\cite{Ciucu}. In Section~\ref{sec:id} this article, we will prove a generalization of Theorem~\ref{thm:yyz} in which we replace $\{a_j\}$ with an arbitrary subset of $A$. The proof of this generalization will follow the spirit of the graphical condensation proofs in~\cite{Kuo}. We can also prove another version of Theorem~\ref{thm:yyz} in which $G$ has an odd number of vertices.

Notice that the relation in Proposition~\ref{prop:int1} can also be rewritten as follows:
\begin{eqnarray}
\lefteqn{M(G) M(G-\{a,b,c,d\}) = M(G-\{a,b\}) M(G-\{c,d\}) + } \nonumber \\
& & M(G-\{a,d\}) M(G-\{b,c\}) - M(G-\{a,c\}) M(G-\{b,d\}). \label{eq:int1pfaff}
\end{eqnarray}
We can interpret the RHS of equation (\ref{eq:int1pfaff}) as the Pfaffian of quantities $M(G-\{i,j\})$ where $i,j \in A=\{a,b,c,d\}$. We can generalize equation (\ref{eq:int1pfaff}) to larger Pfaffians where $A$ is an even set of vertices around a face of $G$. This will be demonstrated in Section~\ref{sec:pfaff}. Finally, in Section~\ref{sec:det}, we will examine some specials cases in which the Pfaffian becomes a determinant.

\section{Additional Condensation Identities}\label{sec:id}

We generalize Theorem~\ref{thm:yyz} as follows:

\begin{thm}\label{thm:2n}
Let $G$ be a plane weighted graph with $2n$ vertices. Let vertices $a_1, b_1, a_2, b_2, \ldots, a_k, b_k (2 \leq k \leq n)$ appear in a cyclic order on a face of $G$, and let $A=\{a_1, a_2, \cdots, a_k\}$ and $B=\{b_1, b_2, \cdots, b_k\}$. Then, no matter how we partition $A$ into disjoint subsets $A_1, A_2$, we have
\begin{equation}\label{eqn:2n}
\sum_{W\subseteq B, |W| \:\mathrm{even}} M(G-W)M(G-A-\overline{W}) 
= \sum_{Y\subseteq B, |Y|-|A_1| \:\mathrm{even}} M(G-A_1-Y)M(G-A_2-\overline{Y}),
\end{equation}
where $\overline{W}=B\backslash W$ and $\overline{Y}=B\backslash Y$, and the second sum ranges over all subsets $Y$ that have the same parity as $A_1$.
\end{thm}

\begin{proof}
Let $\mathcal{H}$ be the set of multigraphs $H$ on the vertices of $G$ that have the following properties:
\begin{enumerate}
\item Each vertex in $A$ and $B$ has degree one, while all other vertices have degree 2. 
\item Each edge in $H$ must also exist in $G$; however, some edges in $H$ may appear twice (thus $\mathcal{H}$ is a set of multigraphs).
\item The length of each cycle in $H$ is even.
\end{enumerate}
The connected components of $H$ include cycles, doubled edges, and paths that run from vertices in $A$ to the vertices in $B$. Since $G$ is a plane graph, these paths cannot intersect, for otherwise some vertex in $H$ will have degree greater than 2. It is also impossible for both endpoints of a path to be in $A$. If that were the case, then that path would divide $H$ into two parts so that the number of vertices in $A$ or $B$ in each part is odd. This contradicts the fact that the number of endpoints of the paths in each part must be even. Similarly, no path can have both endpoints in $B$.

Since $H$ has an even number of vertices, and the cycles and doubled edges occupy an even number of vertices, the number of vertices contained within the paths of $H$ must also be even. Thus the number of even paths (each containing an odd number of vertices) must be even.

Let $k(H)$ be the number of cycles in $H$, and let $w(H)$ be the product of all the edge weights in $H$. For each doubled edge in $H$, the weight is multiplied twice in the product $w(H)$.  We show that both sides of equation (\ref{eqn:2n}) are equal to a third quantity,
\[ S=\sum_{H \in \mathcal{H}} 2^{k(H)} w(H). \]
To show that the LHS of equation (\ref{eqn:2n}) is equal to $S$, we need to show that (1) when we superimpose a matching of $G-W$ onto a matching of $G-A-\overline{W}$, we get a multigraph in $\mathcal{H}$, and (2) each multigraph $H \in \mathcal{H}$ can be partitioned into matchings of $G-W$ and $G-A-\overline{W}$ in $2^{k(H)}$ ways for exactly one even subset $W\subseteq B$.

Consider the resulting graph when we superimpose a matching of $G-W$ onto a matching of $G-A-\overline{W}$. Each vertex in $A$ and $B$ will have degree one, and all other vertices will have degree 2. Also, any cycle must have even length since adjacent edges in the cycle cannot come from the same matching. This resulting multigraph is a member of $\mathcal{H}$. 

We now show that each multigraph $H \in \mathcal{H}$ can be partitioned into matchings of $G-W$ and $G-A-\overline{W}$ in $2^{k(H)}$ ways for exactly one even subset $W\subseteq B$. We would then partition $H$ as follows:
\begin{enumerate}
\item For each doubled edge, include that edge into both matchings.
\item For each path with endpoints $a_i$ and $b_j$, the edge incident to $a_i$ must belong to the matching of $G-W$ since $a_i \not\in G-A-\overline{W}$. Then adjacent edges in the path belong to different matchings.
\item For each cycle, we partition the edges so that adjacent edges belong to different matchings.
\end{enumerate}

If the length of the path between $a_i$ and $b_j$ is odd, then the edge incident to $b_j$ must be in the matching of $G-W$.  Thus the partition is possible if and only if $b_j \not\in W$. On the other hand, if the length of the path between $a_i$ and $b_j$ is even, then the edge incident to $b_j$ must be in the matching of $G-A-\overline{W}$. In that case, the partition is possible if and only if $b_j \in W$. Thus the partition of $H$ is possible if and only if $W$ is the set of vertices $b_i$ that are endpoints of paths of even length. And since there are an even number of even paths in $H$, $W$ is an even subset.

For each cycle, we have two choices for which edges will be put in each matching. Since there are $k(H)$ cycles, there will be $2^{k(H)}$ ways to partition $H$ into matchings of $G-W$ and $G-A-\overline{W}$. Therefore
\[ \sum_{H \in \mathcal{H}} 2^{k(H)} w(H) = \sum_{W\subseteq B, |W| \:\mathrm{even}} M(G-W)M(G-A-\overline{W}).\]

A similar argument can be made to equate $S$ with the RHS of equation (\ref{eqn:2n}). When we superimpose a matching of $G-A_1-Y$ onto a matching of $G-A_2-\overline{Y}$, we also get a multigraph in $\mathcal{H}$. We can also show that any $H \in \mathcal{H}$ can be partitioned into matchings of $G-A_1-Y$ and $G-A_2-\overline{Y}$ in $2^{k(H)}$ ways for exactly one subset $Y\subseteq B$. This time, each edge adjacent to $a_i$ goes into the matching of $G-A_1-Y$ if and only if $a_i \in A_2$. Thus a partition of $H$ is possible if and only if $Y$ is the set of vertices $b_j$ that are either (1) the endpoint of an odd path whose other endpoint is in $A_1$, or (2) the endpoint of an even path whose other endpoint is in $A_2$. The parity of $|Y|$ must be the same as $|A_1|$ so that $G-A_1-Y$ can have perfect matchings. Therefore
\[ \sum_{H \in \mathcal{H}} 2^{k(H)} w(H) = \sum_{Y\subseteq B, |Y|-|A_1| \:\mathrm{even}} M(G-A_1-Y)M(G-A_2-\overline{Y}), \]
and so equation (\ref{eqn:2n}) is proved.
\end{proof}

\begin{rem}
If we set $A_1=\emptyset$ and $A_2=A$, equation (\ref{eqn:2n}) becomes reflexive.
\end{rem}

\begin{rem}
If we set $A_1$ to contain a single vertex $a_j$, then we derive Theorem~\ref{thm:yyz}.
\end{rem}

By letting $G$ be bipartite in Theorem~\ref{thm:2n}, we derive the following corollary:

\begin{cor}
Let $G=(U,V)$ be a plane weighted bipartite graph in which $U=\{u_i|1\leq 1 \leq n\}$ and $V=\{v_i|1\leq 1 \leq n\}$. Let vertices $a_1, b_1, a_2, b_2, \ldots, a_k, b_k (2 \leq k \leq n)$ appear in a cyclic order on a face of $G$, and let $A=\{a_1, a_2, \cdots, a_k\}\subseteq U$ and $B=\{b_1, b_2, \cdots, b_k\} \subseteq V$. Then no matter how we partition $A$ into disjoint subsets $A_1, A_2$, we have
\begin{equation}
M(G)M(G-A-B) = \sum_{Y\subseteq B, |Y|=|A_1|} M(G-A_1-Y)M(G-A_2-\overline{Y}),
\end{equation}
where $\overline{Y}=B\backslash Y$, and the second sum ranges over all subsets $Y\subseteq B$ with the same cardinality as $A_1$.
\end{cor}

\begin{proof}
This identity is a special case of Theorem~\ref{thm:2n}. For any nonempty subset $W \subseteq B$, $M(G-W)=0$ since $G-W$ has unequal numbers of vertices in $U$ and $V$. Similarly, $M(G-A_1-Y)=0$ whenever $|Y| \neq |A_1|$.
\end{proof}

\begin{rem}
If we set $A_1$ to contain a single vertex $a_j$ and $Y$ range over single vertices in $B$, then we derive Corollary 2.3 of~\cite{YYZ}.
\end{rem}

The following theorem is analogous to Theorem~\ref{thm:2n} in which $G$ has an odd number of vertices:

\begin{thm}\label{thm:2n+1}
Let $G$ be a plane weighted graph with $2n+1$ vertices. Let vertices $a_1, b_1, a_2, b_2, \ldots, a_k, b_k (2 \leq k \leq n)$ appear in a cyclic order on a face of $G$, and let $A=\{a_1, a_2, \cdots, a_k\}$ and $B=\{b_1, b_2, \cdots, b_k\}$. Then, no matter how we partition $A$ into disjoint subsets $A_1, A_2$, we have
\begin{equation}
\sum_{W\subseteq B, |W| \:\mathrm{odd}} M(G-W)M(G-A-\overline{W}) 
= \sum_{Y\subseteq B, |Y|-|A_1| \:\mathrm{odd}} M(G-A_1-Y)M(G-A_2-\overline{Y}), 
\end{equation}
where $\overline{W}=B\backslash W$ and $\overline{Y}=B\backslash Y$, and the second sum ranges over all subsets $Y$ whose parity is the opposite of $A_1$.
\end{thm}

\begin{proof}
The proof basically follows that of Theorem~\ref{thm:2n}. We define $\mathcal{H}$ as we did in the proof of Theorem~\ref{thm:2n}. The only property that a multigraph $H \in \mathcal{H}$ differs is that since $H$ has an odd number of vertices, $H$ also has an odd number of even paths. Then any $H$ can be partitioned into matchings of $G-W$ and $G-A-\overline{W}$ for exactly one odd set $W$. We could also partition $H$ into matchings of $G-A_1-Y$ and $G-A_2-\overline{Y}$ for exactly one set $Y$ for which $|Y|$ and $|A_1|$ have opposite parity. Each partition can be done in $2^{k(H)}$ ways, where $k(H)$ is the number of cycles in $H$.
\end{proof}

The following corollary is an adaptation of Proposition~\ref{prop:int1} for when $G$ has an odd number of vertices.

\begin{cor}
Let $G$ be a plane weighted graph with $2n+1$ vertices. Let vertices $a_1, b_1, a_2, b_2$ appear in a cyclic order on a face of $G$. Then
\begin{eqnarray}
\lefteqn{M(G-a_1)M(G-\{a_2,b_1,b_2\}) + M(G-a_2)M(G-\{a_1,b_1,b_2\}) } \nonumber \\
& & = M(G-b_1)M(G-\{a_1,a_2,b_2\}) + M(G-b_2)M(G-\{a_1,a_2,b_1\}).
\end{eqnarray}
\end{cor}

\begin{rem}
If $G=(U,V)$ is a bipartite graph in which $|U|=|V|+1$, and $a_1,b_1,a_2 \in U$ and $b_2 \in V$, then we derive Theorem 2.4 of ~\cite{Kuo}:
\begin{equation}
M(G-a_1)M(G-\{a_2,b_1,b_2\}) + M(G-a_2)M(G-\{a_1,b_1,b_2\}) = M(G-b_1)M(G-\{a_1,a_2,b_2\}).
\end{equation}
\end{rem}

This next corollary is a combination of special cases of Theorems~\ref{thm:2n} and~\ref{thm:2n+1}:

\begin{cor}\label{cor:n+k}
Let $G=(U,V)$ be a plane weighted bipartite graph in which $U=\{u_i|1\leq 1 \leq n+k\}$ and $V=\{v_i|1\leq 1 \leq n\}$. Let vertices $a_1, b_1, a_2, b_2, \ldots, a_k, b_k (2 \leq k \leq n)$ appear in a cyclic order on a face of $G$, and let $A=\{a_1, a_2, \cdots, a_k\}\subseteq U$ and $B=\{b_1, b_2, \cdots, b_k\}\subseteq U$. Then no matter how we partition $A$ into disjoint subsets $A_1, A_2$, we have
\begin{equation}
M(G-A)M(G-B) = \sum_{Y\subseteq B, |Y|=k-|A_1|} M(G-A_1-Y)M(G-A_2-\overline{Y}), 
\end{equation}
where $\overline{Y}=B\backslash Y$, and the second sum ranges over all subsets $Y$ such that $|Y|+|A_1|=k$.
\end{cor}

\begin{proof}
As another special case of Theorem~\ref{thm:2n}, the only subset $W\subseteq B$ for which $G-W$ has perfect matchings is $W=B$. Similarly, $M(G-A_1-Y)=0$ whenever $|A_1|+|Y| \neq k$.
\end{proof}

\begin{rem}
If we set $k=2$ in Corollary~\ref{cor:n+k}, we derive Theorem 2.5 in~\cite{Kuo}.
\end{rem}

\section{Pfaffian Generalization of Graphical Condensation}\label{sec:pfaff}

In order to generalize equation (\ref{eq:int1pfaff}) to larger Pfaffians, we will need to set up a series of hypotheses for our generalization. Let $G$ be a plane graph. For an even number $n$, let $A$ be a set of vertices $a_1, a_2, \ldots, a_n$ that appear in a cyclic order on a face of $G$. Further, divide $A$ into two disjoint subsets $A_K$ and $A_H$ and let $H = G-A_K$. For each pair $(i,j)$, $1 \leq i < j \leq n$, define $H_{ij}$ to be one of the following subgraphs of $G$:
\begin{itemize}
\item If $a_i, a_j \in A_K$, let $H_{ij}$ be induced by $H$ with $a_i$ and $a_j$.
\item If $a_i \in A_K$ and $a_j \in A_H$, let $H_{ij}$ be induced by $H-a_j$ and $a_i$. 
\item If $a_i \in A_H$ and $a_j \in A_K$, let $H_{ij}$ be induced by $H-a_i$ and $a_j$. 
\item If $a_i, a_j \in A_H$, let $H_{ij} = H-a_i-a_j$.
\end{itemize}
In particular, note that $a_i \in H_{ij}$ if and only if $a_i \not\in H$.

\begin{thm}\label{thm:pfaff}
If $H$ has exactly one matching (i.e. $M(H)=1$), then the number of perfect matchings of $K = G-A_H$ is
\begin{equation}\label{eq:pfaff}
M(K) = M(G-A_H) = \Pf(M(H_{ij}))_{1 \leq i \leq j \leq n}.
\end{equation}
\end{thm}

\begin{proof}
Let $M_H$ be the lone perfect matching of $H$. Define an \emph{$M_H$-alternating path} from $a_i$ to $a_j$ to be a path in $G$ in which every other edge is in $M_H$, and if either $a_i$ or $a_j \in A_H$, then the path must include the edge from $M_H$ incident with that vertex.

\begin{lem}\label{lem:hij}
The number of $M_H$-alternating paths from $a_i$ to $a_j$ is equal to $M(H_{ij})$.
\end{lem}

\begin{proof}
We establish a bijection between $M_H$-alternating paths from $a_i$ to $a_j$ and matchings of $H_{ij}$. Start with a matching of $H_{ij}$, and consider what happens when we superimpose $M_H$ on that matching of $H_{ij}$. (If either $a_i$ or $a_j$ is not in $H_{ij}$, include it in the combined graph.) A path arises from $a_i$ to $a_j$ since $a_i$ and $a_j$ are the only vertices of degree one. Every other edge in the path is from $M_H$. 
Thus we can associate each perfect matching of $H_{ij}$ with an $M_H$-alternating path from $a_i$ to $a_j$.

Given an $M_H$-alternating path from $a_i$ to $a_j$, we can recover the corresponding matching of $H_{ij}$. First remove all the edges in the path that are from $M_H$, deleting $a_i$ or $a_j$ if necessary. The remaining vertices off the path can be matched in only one way, namely with the edges in $M_H$ that are not on the path. If there were another way to match those remaining vertices, $M_H$ would not be the only perfect matching of $H$. The bijection is established, so the number of $M_H$-alternating paths from $a_i$ to $a_j$ is $M(H_{ij})$.
\end{proof}

We will refer to a collection of $M_H$-alternating paths among the vertices in $A$ as a \emph{nest.}

\begin{lem}\label{lem:mk}
The number of nests of non-intersecting $M_H$-alternating paths through $G$ with the vertices in $A$ as endpoints is equal to $M(K) = M(G-A_H)$.
\end{lem}

\begin{proof}
Just as in the previous lemma, we establish a bijection between matchings of $K=G-A_H$ with nests of non-intersecting $M_H$-alternating paths among $A$.
Consider what happens when we superimpose $M_H$ on a matching $M_K$ of $K$. In the resulting graph $M_K \cup M_H$, the degree of each vertex in $A$ is 1, and the degree is 2 for all other vertices. Therefore there must be $n$ paths with endpoints in $A$. Each path $P$ alternates between edges in $M_H$ and $M_K$. Paths cannot intersect in $M_K \cup M_H$ since some vertex would have degree greater than 2. Thus we can associate each perfect matching of $K$ with a nest of non-intersecting $M_H$-alternating paths among vertices in $A$.

Now start with a nest of non-intersecting $M_H$-alternating paths through $G$ with endpoints in $A$. We will create a perfect matching in $K$ by first removing all edges in $M_H$ from the paths, removing vertices in $A$ if necessary. The remaining vertices off the path can be matched in only one way, namely with the edges in $M_H$ that are not on the path. If there were another way to match those remaining vertices, $M_H$ would not be the only perfect matching of $H$. The bijection is established, so the number of nests of non-intersecting $M_H$-alternating paths is $M(K)$.
\end{proof}

\emph{Proof of Theorem~\ref{thm:pfaff} continued}: Now we establish the Pfaffian identity (\ref{eq:pfaff}). By the definition of a Pfaffian
\begin{equation}\label{eq:pf2}
\Pf(M(H_{ij}))_{1 \leq i \leq j \leq n} = \sum_{F \in \mathcal{F}_n} (-1)^{\chi(F)} \prod_{(i,j)\in F} M(H_{i\sigma(i)}).
\end{equation}
The set $\mathcal{F}_n$ consists of all partitions of $\{1, 2, \ldots, n\}$ into pairs of elements. Each partition $F \in \mathcal{F}_n$ is called a \emph{1-factor}. The \emph{crossing number} $\chi(F)$ is the number of pairs $(a,b), (c,d)$ in $F$ for which $a<c<b<d$. If we place the numbers $1, 2, \ldots, n$ around a circle and draw a chord for each pair in $F$, then $\chi(F)$ is the number of pairs of intersecting chords in the circle.

For each 1-factor $F \in \mathcal{F}_n$, the product in the RHS of equation (\ref{eq:pf2}) represents the number of nests of $M_H$-alternating paths from $a_i$ to $a_j$ for each pair $(i,j) \in F$. Depending on whether the crossing number $\chi(F)$ of the 1-factor is odd or even, each nest will contribute either +1 or $-1$ to the sum. Each nest of non-intersecting paths will contribute +1 to the sum since the crossing number of the corresponding 1-factor is 0. We show how the contributions from the intersecting nests cancel each other out by establishing a bijection between intersecting nests that contribute $-1$ with intersecting nests that contribute +1. To do so, we need to prove a property about intersecting $M_H$-alternating paths.

\begin{prop}\label{prop:ijkl}
Let $\mathcal{C}$ be a nest of intersecting $M_H$-alternating paths $\{P_{ij}: (i,j) \in F\}$ in which $P_{ij}$ is a path from $a_i$ to $a_j$. Suppose paths $P_{ij}$ and $P_{k\ell}$ intersect in $\mathcal{C}$. Let $Q$ be the union of all the edges in $P_{ij}$ and $P_{k\ell}$, keeping all doubled edges where they occur. We show that $Q$ can be partitioned either into paths $P'_{ik}$ and $P'_{j\ell}$ or paths $P'_{i\ell}$ and $P'_{jk}$, but not both. In particular, if $i < j < k < \ell$, then the paths must be $P'_{ik}$ and $P'_{j\ell}$.
\end{prop}

\begin{proof}
Create a simple graph $R$ in which $V(R) = V(Q)$ and each edge is from $P_{ij}$ or $P_{k\ell}$. Color an edge red if it belongs only to $P_{ij}$, blue if it belongs only to $P_{k\ell}$, and purple if it belongs to both paths. Only vertices $a_i, a_j, a_k$, and $a_{\ell}$ have degree one in $R$. At each vertex of degree two, both incident edges have the same color. When a vertex $v$ has a degree greater than two, it must be an intersection of both paths. One edge $e$ incident to $v$ must be in $M_H$, and it must be shared by both paths; thus $e$ is a purple edge. At most two other edges can be incident to $v$, one red and one blue. Thus the maximum degree of any vertex is three, and one red, one blue, and one purple $M_H$ edge meet at those vertices. We will define a \emph{section} of $R$ to be a connected component of edges of the same color. Note that each section is an $M_H$-alternating path between two vertices of degree 1 or 3. Both ends of a purple section are at vertices of degree 3 with edges from $M_H$. If a red or blue section ends at a vertex of degree 3, then the final edge is not in $M_H$.

If you delete all the purple edges from $R$, you are left with several connected components of red and blue edges. In particular, two of these components are paths between $a_i$, $a_j$, $a_k$, and $a_{\ell}$. Call these paths \emph{red-blue paths}. Since they are separate components, these red-blue paths do not intersect. Neither red-blue path is purely red or blue; if that were the case, those red-blue paths would become all of $P_{ij}$ and $P_{k\ell}$, which intersect.

\begin{lem}\label{lem:ends}
It is impossible for $a_i$ and $a_j$ to be the endpoints of one red-blue path, and $a_k$ and $a_{\ell}$ to be endpoints of the other.
\end{lem}

\begin{proof}
Assume the contrary, that $a_i$ and $a_j$ are endpoints of one red-blue path and $a_k$ and $a_{\ell}$ are endpoints of the other. Then the first and last edges of one red-blue path are red, and the first and last edges of the other red-blue path are blue. Now swap the colors in the red-blue path with $a_k$ and $a_{\ell}$. Then delete all the red edges from this recolored graph, along with $a_i$, $a_j$, $a_k$, and $a_{\ell}$. Since neither red-blue path was purely red or blue, there must be some blue edges remaining from those red-blue paths. Since there are no more vertices of degree one or three, the remaining blue and purple edges must form at least one cycle $C$. This cycle is composed of blue and purple sections, each of which is $M_H$-alternating. At a vertex where a blue and a purple edge meet, the purple edge is in $M_H$ and the blue one is not. Thus the entire cycle is $M_H$-alternating. But by swapping membership in $M_H$ among the edges of $C$, we create another perfect matching for $H$! This contradicts our assumption that $H$ has only one perfect matching, so the lemma follows.
\end{proof}

As a corollary to Lemma~\ref{lem:ends}, the first and last edges of each red-blue path must be different colors. Assume without loss of generality that $a_i$ and $a_{\ell}$ are the endpoints of one red-blue path, and $a_j$ and $a_k$ are the endpoints of the other. We now show that $Q$ cannot be partitioned into $M_H$-alternating paths $P'_{i\ell}$ and $P'_{jk}$ from $a_i$ to $a_{\ell}$ and from $a_j$ to $a_k$. Suppose that such a partition is possible. At every junction of red, blue, and purple sections in $R$, the purple section is shared by both paths, and the red and blue sections belong to different paths. The red section beginning with $a_i$ in $R$ must belong to $P'_{i\ell}$. So as we trace the red-blue path from $a_i$, each red section must belong to $P'_{i\ell}$ and each blue section must belong to $P'_{jk}$. However, the final section of the red-blue path is both blue and ends at $a_{\ell}$. Thus we have a blue section also belonging to $P'_{i\ell}$, which is impossible. Thus the partition into $P'_{i\ell}$ and $P'_{jk}$ is impossible.

On the other hand, we show that $Q$ can indeed be partitioned into paths $P'_{ik}$ and $P'_{j\ell}$ from $a_i$ to $a_k$ and from $a_j$ to $a_{\ell}$. To do this, we swap the color of each edge in the red-blue path from $a_j$ to $a_k$ in $R$ and call the new coloring $R'$. In $R'$ each vertex of degree 2 is still incident to two edges of the same color, and each vertex of degree 3 is still incident to a red, blue, and purple edge.

\begin{lem}\label{lem:re-part}
If we delete all the blue edges in $R'$, the remaining red and purple edges form a single $M_H$-alternating path $P'_{ik}$ from $a_i$ to $a_k$.
\end{lem}

\begin{proof}
Except for $a_i$ and $a_k$, each vertex will be connected to either zero or two edges. Therefore if we start tracing a path from $a_i$ along the remaining edges, we must end at $a_k$. This path is $M_H$-alternating since each red and purple section is $M_H$-alternating, and at each vertex incident to a red and purple edge, the purple edge is in $M_H$ and the red one is not. Any remaining edges not on this path must form cycles. But as in the proof of Lemma~\ref{lem:ends}, this cycle is also $M_H$-alternating, meaning that another matching of $H$ exists. This is a contradiction, so the path $P'_{ik}$ contains all red and purple edges in $R'$.
\end{proof}

\emph{Proof of Proposition~\ref{prop:ijkl} continued}: Likewise, if we delete all the red edges in $R'$, the remaining blue and purple edges form a single $M_H$-alternating path $P'_{j\ell}$ from $a_j$ to $a_{\ell}$. Together, paths $P'_{ik}$ and $P'_{j\ell}$ cover each red and blue edge once and each purple edge twice. These paths form a partition of $Q$.

Finally, assume that $i<j<k<\ell$ and that paths $P_{ij}$ and $P_{k\ell}$ intersect in $\mathcal{C}$. Then from Lemma~\ref{lem:ends} and the fact that red-blue paths in $R$ do not intersect, the red-blue paths must run from $a_i$ to $a_{\ell}$ and from $a_j$ to $a_k$. Our argument showed that we can partition $Q$ only into paths $P'_{ik}$ and $P'_{j\ell}$.
\end{proof} 

\emph{Proof of Theorem~\ref{thm:pfaff} continued}: Given a nest $\mathcal{C}$ of intersecting $M_H$-alternating paths $\{P_{ij}: (i,j) \in F\}$, let $i$ be the smallest integer for which the path $P_{ij}$ intersects another path. As we start from $a_i$ and trace $P_{ij}$, let $v$ be the first vertex to be shared with another path. Let $\ell$ be the smallest integer for which $P_{k\ell}$ also goes through $v$, and a red-blue path exists from $a_i$ to $a_{\ell}$. Then the union of $P_{ij}$ and $P_{k\ell}$ can be partitioned into paths $P'_{ik}$ and $P'_{j\ell}$ (but not into $P'_{i\ell}$ and $P'_{jk}$). Paths $P'_{ik}$ and $P'_{j\ell}$ are created using the procedure described in Lemma~\ref{lem:re-part}. Next, we create another nest $\mathcal{C}'$ of paths which include $P'_{ik}$, $P'_{j\ell}$, and all paths in $\mathcal{C}$ except $P_{ij}$ and $P_{k\ell}$. Let $\mathcal{C}''$ be the result of applying the same transformation to this new nest $\mathcal{C}'$. The smallest integer $i'$ for which $a_{i'}$ lies on an intersecting path is still $i'=i$. Vertex $v$ is still the first vertex along $P'_{ik}$ to be shared by another path, and $\ell'=\ell$ is still the smallest integer for which $P_{j\ell}$ also goes through $v$, and a red-blue path exists from $a_i$ to $a_{\ell}$. Using the procedure from Lemma~\ref{lem:re-part} once again, we receive the old paths $P_{ij}$ and $P_{k\ell}$, so $\mathcal{C}'' = \mathcal{C}$.

To complete the proof, we examine the crossing numbers of $F$ and $F'$, where $F'$ is the same as $F$ except that it includes pairs $(i,k)$ and $(j,\ell)$ but not $(i,j)$ and $(k,\ell)$. Assume WLOG that $i<j<k<\ell$. If we use the chord-crossing model for $F$, consider how the chord $(u,v) \in F$ separates $i$, $j$, $k$, and $\ell$ into two sides and crosses the chords connecting them in $F$ and $F'$. If all four points are on the same side, then $(u,v)$ does not either chord in $F$ or $F'$. If $(u,v)$ separates one point from the other three, then $(u,v)$ crosses exactly one of the chords in both $F$ and $F'$. If $(u,v)$ separates $i$ and $j$ from $k$ and $\ell$, then $(u,v)$ crosses neither chord in $F$ and both chords in $F'$. Finally, if $(u,v)$ separates $i$ and $\ell$ from $j$ and $k$, then $(u,v)$ crosses both chords in $F$ and $F'$. Thus the parity of crossings with $(u,v)$ remains unchanged. However, $(i,j)$ and $(k,\ell)$ do not cross in $F$ but $(i,k)$ and $(j,\ell)$ do in $F'$. Therefore
\[ (-1)^{\chi(F)} = -1 \cdot (-1)^{\chi(F')}, \]
and so the contributions from $\mathcal{C}$ and $\mathcal{C}'$ cancel each other out in the sum in equation (\ref{eq:pf2}). The only remaining contributions are from the $M(G)$ nests of non-intersecting paths, and relation (\ref{eq:pfaff}) is proved. 
\end{proof} 

The following corollaries are the special cases of Theorem~\ref{thm:pfaff} in which $A=A_K$ or $A=A_H$.

\begin{cor}
Let $a_1, a_2, \ldots, a_n$ be an even number of vertices that appear in a cyclic order along a face of a plane graph $G$. Let $A=\{a_1, a_2, \cdots, a_n\}$ and $H = G-A$. Define $H_{ij}$ to be the subgraph induced by the vertices in $H$ along with $a_i$ and $a_j$. If $H$ has exactly one perfect matching, then 
\begin{equation}
M(G) = \Pf(M(H_{ij}))_{1 \leq i \leq j \leq n}.
\end{equation}
\end{cor}

\begin{cor}
Let $a_1, a_2, \ldots, a_n$ be an even number of vertices that appear in a cyclic order along a face of a plane graph $G$. Let $A=\{a_1, a_2, \cdots, a_n\}$ and $K = G-A$. Define $H_{ij} = G-a_i-a_j$. If $G$ has exactly one perfect matching, then 
\begin{equation}
M(K) = M(G-A) = \Pf(M(H_{ij}))_{1 \leq i \leq j \leq n}.
\end{equation}
\end{cor}

\section{Determinant Generalization of Graphical Condensation}\label{sec:det}

In this section, we examine some corollaries to Theorem~\ref{thm:pfaff} in which $G=(U,V)$ is a plane bipartite graph. This version expresses the number of perfect matchings of $G$ as a determinant whose entries are numbers of perfect matchings of subgraphs of $G$. The following theorem bears a strong resemblance to a theorem by Gessel and Viennot~\cite{GV} in which the number of non-intersecting lattice paths is expressed as a determinant formula. 

Let $G=(U,V)$ be a plane bipartite graph, and vertices $a_1, a_2, \ldots, a_n, b_n, b_{n-1}, \ldots, b_2, b_1$ appear in a cyclic order on a face of $G$. Let $A = \{a_1, a_2, \cdots, a_n\}$ and $B=\{b_1, b_2, \cdots, b_n\}$, and let $L = G-(A \cap U)-(B \cap V)$. This time, for each pair $(i,j)$, $1 \leq i,j \leq n$, we define $L_{ij}$ as follows:
\begin{itemize}
\item If $a_i, b_j \in U$, let $L_{ij}$ be induced by $L-a_j$ and $a_i$.
\item If $a_i \in U$ and $b_j \in V$, let $L_{ij}$ be induced by $L$ with $a_i$ and $a_j$.
\item If $a_i \in V$ and $b_j \in U$, let $L_{ij} = L-a_i-a_j$.
\item If $a_i, b_j \in V$, let $L_{ij}$ be induced by $L-a_i$ and $a_j$.
\end{itemize}

\begin{thm}
Let $K = G-(A \cap V)-(B \cap U)$. If $L$ has exactly one perfect matching, then 
\begin{equation}\label{eq:det}
M(K) = \det[M(L_{ij})]_1^n.
\end{equation}
\end{thm}

\begin{proof}
We apply Theorem~\ref{thm:pfaff} in which we let $a_{n+j} = b_{n+1-j}$ for $1 \leq j \leq n$, and let $A_H = (A \cap V) \cup (B \cap U)$ and $A_K = (A \cap U) \cup (B \cap V)$. If $i \leq n$, then $a_i \in U$ if and only if $a_i \in A_K$. Also, if $i>n$, then $a_i \in U$ if and only if $a_i \in A_H$. Therefore if both $i,j \leq n$, then $H_{ij}$ has two more vertices in $U$ than $V$. Similarly, if both $i,j > n$, then $H_{ij}$ has two more vertices in $V$ than $U$. Thus $M(H_{ij})=0$ when both $i,j \leq n$, or both $i,j >n$. When those quantities are set to 0 in a Pfaffian, the remaining terms form a determinant. Specifically,
\begin{equation}
M(G) = \Pf(M(H_{ij}))_{1\leq i<j \leq n} = \det[M(H_{i(2n+1-j)}]_1^n = \det[M(L_{ij})]_1^n.
\end{equation}
\end{proof}

Alternatively, we could adapt the proof of Theorem~\ref{thm:pfaff} for this special case.
Let $M_L$ be the lone perfect matching of $L$. As in Lemmas~\ref{lem:hij} and~\ref{lem:mk}, $M(L_{ij})$ is also the number of $M_L$-alternating paths in $G$ from $a_i$ to $b_j$, and $M(G)$ is also the number of non-intersecting nests of paths in $G$ with endpoints in $A \cup B$. Since it is impossible for the endpoints of an $M_L$-alternating path to be both in $A$ or both in $B$, each path must have one endpoint in $A$ and the other in $B$. Let $P_i$ denote the $M_L$-alternating path starting at $a_i$. We now demonstrate that $P_i$ must end at $b_i$ for each $i=1, \ldots, n$. For suppose that $j$ is the smallest integer for which the path $P_j$ does not end at $b_j$. Let the paths starting at $a_j$ and $b_j$ end at $b_k$ and $a_{\ell}$, respectively, where $j<k,\ell$. Since vertices $a_j, a_{\ell}, b_k$, and $b_j$ appear in that order on a face of $G$, it is impossible for paths to exist from $a_j$ to $b_k$ and from $a_{\ell}$ to $b_j$ without intersecting. Thus each path $P_i$ must end at $b_i$. We can now associate each perfect matching of $G$ with a nest of $n$ non-intersecting $M_L$-alternating paths from $a_i$ to $b_i$ for $i=1, \ldots, n$.

Now we establish the determinant identity (\ref{eq:det}). Let $S_n$ be the symmetric group on $n$ elements, and let $\sgn(\sigma)$ be the signature of the permutation $\sigma \in S_n$. Then by the definition of a determinant,
\begin{equation}\label{eq:det2}
\det [M(L_{ij})]_1^n = \sum_{\sigma \in S_n} \sgn(\sigma) \prod_{i=1}^n M(L_{i\sigma(i)}).
\end{equation}
For each permutation $\sigma \in S_n$, the product in the RHS of equation (\ref{eq:det2}) represents the number of nests of $M_L$-alternating paths from $a_i$ to $b_{\sigma(i)}$. Depending on the signature of the permutation, each nest will contribute either +1 or $-1$ to the sum. Each nest of non-intersecting paths will contribute +1 to the sum since the signature of the identity permutation is +1. We show how the contributions from the intersecting nests cancel each other out. We establish a bijection between  intersecting nests that contributes $-1$ with intersecting nests that contribute +1. Let $\mathcal{C}$ be a nest of intersecting paths $\{P_1, \ldots, P_n\}$ in which $P_i$ is a path from $a_i$ to $b_{\sigma(i)}$. Let $i$ be the smallest integer for which the path $P_i$ starting at $a_i$ intersects another path. As we start from $a_i$ and trace $P_i$, let $v$ be the first vertex to be shared with another path. Let $j$ be the smallest integer greater than $i$ for which $P_j$ also goes through $v$. Next, we create another nest $\mathcal{C}'$ of paths $\{P'_1,\ldots, P'_n\}$ where $P'_k=P_k$ for all $k \neq i,j$. In $\mathcal{C}'$, path $P'_i$ follows $P_i$ from $a_i$ to $v$, and $P_j$ from $v$ to $b_{\sigma(j)}$. Path $P'_j$ follows $P_j$ from $a_j$ to $v$, and $P_i$ from $v$ to $b_{\sigma(i)}$. The corresponding permutation $\sigma'$ for $\mathcal{C}'$ is the same as $\sigma$ except $\sigma'(i) = \sigma(j)$ and $\sigma'(j) = \sigma(i)$. Because we have swapped two places in $\sigma$ to create $\sigma'$, we have $\sgn(\sigma') = -\sgn(\sigma)$. Our bijection maps $\mathcal{C}$ to $\mathcal{C}'$, and their contributions to the sum in equation (\ref{eq:det2}) cancel out. The only remaining contributions are from the $M(G)$ nests of non-intersecting paths, and relation (\ref{eq:det}) is proved.

Two corollaries result if $A \subseteq U$ and $B \subseteq V$, or vice versa.

\begin{cor}\label{cor:det1}
Let $G=(U,V)$ be a plane bipartite graph, and vertices $a_1, a_2, \ldots, a_n, b_n, b_{n-1}, \ldots, b_2, b_1$ appear in a cyclic order on a face of $G$. Let $A=\{a_1, a_2, \cdots, a_n\}\subseteq U$ and $B=\{b_1, b_2, \cdots, b_n\}\subseteq V$, and let $L = G-A-B$. Define $L_{ij}$ to be the subgraph induced by the vertices in $L$ along with $a_i$ and $b_j$. If $L$ has exactly one perfect matching, then 
\begin{equation}
M(G) = \det[M(L_{ij})]_1^n.
\end{equation}
\end{cor}

\begin{cor}\label{cor:det2}
Let $G=(U,V)$ be a plane bipartite graph, and vertices $a_1, a_2, \ldots, a_n, b_n, b_{n-1}, \ldots, b_2, b_1$ appear in a cyclic order on a face of $G$. Let $A=\{a_1, a_2, \cdots, a_n\} \subseteq V$ and $B=\{b_1, b_2, \cdots, b_n\} \subseteq U$, and let $K = G-A-B$. Define $L_{ij} = G-a_i-b_j.$ If $G$ has exactly one perfect matching, then 
\begin{equation}
M(K) = \det[M(L_{ij})]_1^n.
\end{equation}
\end{cor}


\begin{thebibliography}{99}
\bibitem{Carlitz}
L. Carlitz, Rectangular arrays and plane partitions,
{\it Acta Arithmetica} {\bf 13} (1967), 29--47.

\bibitem{Ciucu}
M. Ciucu, Enumeration of perfect matchings in graphs with reflective symmetry, J. Combin. Theory Ser. A, 77 (1997), 67-97

\bibitem{GV} 
I. Gessel, G. Viennot. Binomial determinants, paths, and hook length
formulae, {\it Advances in Mathematics} {\bf 58} (1985), no 3: 300--321.

\bibitem{Kuo} 
E. Kuo, Applications of Graphical Condensation for Enumerating Matchings and Tilings, Theoret. Comput. Sci., 319 (2004), 29-57.


\bibitem{Propp} 
J. Propp, Generalized Domino-Shuffling, Theoret. Comput. Sci., 303 (2003), 267-301.

\bibitem{YYZ} 
W.G. Yan, Y.N. Yeh, F.J. Zhang, Graphical condensation of plane graphs: a combinatorial approach, Theoret. Comput. Sci., 349 (2005), 452-461.

\bibitem{YZ} 
W.G. Yan, F.J. Zhang. Graphical condensation for enumerating perfect matchings, J. Combin. Theory Ser. A, 110 (2005), 113-125.
\end{thebibliography}

\end{document}